\documentclass{article}
\usepackage{amsfonts}
\usepackage{amsmath}
\usepackage{amssymb}

\newtheorem{theorem}{Theorem}[section]
\newtheorem{lemma}[theorem]{Lemma}
\newtheorem{proposition}[theorem]{Proposition}
\newtheorem{corollary}[theorem]{Corollary}

\newtheorem{example}[theorem]{Example}
\newtheorem{remark}[theorem]{Remark}

\newtheorem{corolario}[theorem]{Corollary}
\newenvironment{proof}[1][Proof]{\noindent\textbf{#1.}}{\ $\square$\bigskip}
\begin{document}

\title{Global decomposition of a Lorentzian manifold as a Generalized
Robertson-Walker space\thanks{%
This paper was supported in part by MEYC-FEDER Grant MTM2007-60016.}}
\author{Manuel Guti\'errez \\
Departamento de \'{A}lgebra, Geometr\'{\i}a y Topolog\'{\i}a.\\
Universidad de M\'alaga. Spain. mgl@agt.cie.uma.es \and Benjam\'in Olea \\
Departamento de \'{A}lgebra, Geometr\'{\i}a y Topolog\'{\i}a.\\
Universidad de M\'alaga. Spain. benji@agt.cie.uma.es}
\date{}
\maketitle

\begin{abstract}
Generalized Robertson-Walker (GRW) spaces constitute a quite important family in Lorentzian geometry, and it is an interesting question to know
whether a Lorentzian manifold can be decomposed in such a way. It is well known that the existence of a suitable vector field guaranties the
local decomposition of the manifold. In this paper, we give conditions on the curvature which ensure a global decomposition and apply them to
several situations where local decomposition appears naturally. We also study the uniqueness question, obtaining that the de Sitter spaces are
the only non trivial complete Lorentzian manifolds with more than one GRW decomposition. Moreover, we show that the Friedmann Cosmological
Models admit an unique GRW decomposition, even locally.
\end{abstract}



\section{Introduction}

If $I\subset \mathbb{R}$ is an open interval, $(L,g_{0}$) a Riemannian manifold and $f\in C^{\infty }(I)$, the warped product $I\times _{f}L$ is
the product manifold $I\times L$ with metric $-dt^2+f(t)^2g_{0}$ and is called a Generalized Robertson-Walker space, GRW space in short. These
spaces were introduced in \cite{alias} and they have been widely studied since then. It is well known that if a Lorentzian manifold admits a
timelike, closed and conformal vector field, then it is locally a GRW space, \cite{Harris2,Sanchez}, but the presence of such a vector field
does not allow us to decide if the decomposition is global. The best we can say, under suitable completeness conditions, is that the manifold is
a quotient
of a GRW space, which is not always a global product, see examples \ref%
{ejemplo1} and \ref{ejemplo2}.

In this paper we study purely geometric conditions to achieve a global decomposition as a GRW space and apply them to obtain global GRW from
local ones. For example, any Robertson-Walker space is a perfect fluid, so it is natural to find out under which conditions a perfect fluid is a
Robertson-Walker space. Other situation is suggested by photon surfaces in General Relativity, which are timelike totally umbilic submanifolds,
\cite{Ellis}. In a local GRW space, photon surfaces inherit the local GRW structure of the ambient space. We study how far they are in fact
global GRW spaces.

Finally, we investigate the GRW decomposition uniqueness obtaining that the
de Sitter spaces are the unique complete Lorentzian manifolds with more than
one (with non constant warping function) GRW decomposition. On the other
hand, in a non necessarily complete manifold, if the ligthlike sectional
curvature is not zero for any degenerate planes at a point $p$, then there
is at most one local GRW structure in a neighborhood of $p$.

The non complete case is interesting in Cosmology. In fact, Friedmann spaces are incomplete GRW spaces with a distinguished family of comoving
observers which physically represents the average galaxy evolution in the spacetime. From a mathematical point of view, comoving observers are
integral curves of the unitary of a timelike, closed and conformal vector field which gives the global decomposition as a GRW space. Since the
ligthlike sectional curvature of the Friedmann spaces is always positive, this decomposition is unique, even locally, which means that the
Friedmann mathematical representation of the galaxy evolution is unambiguous.

\section{Preliminaries}

As we said, we are going to study when a space splits as a GRW space $M=I\times_f L$. In this space there are two distinguished vector fields.
The first one is $U=\frac{\partial}{\partial t}$ which is a \textit{reference frame} i. e., a timelike and unit vector field. Moreover, it is
closed (its metrically equivalent one form is closed) and orthogonally conformal ($L_Ug(X,Y)=2\alpha g(X,Y)$ for all $X,Y\in U^\perp $ for some
function $\alpha\in C^\infty(M)$). The other one is $V=f\frac{\partial}{\partial t}$ which is timelike, closed and conformal. Our interest is to
study the implications that the existence of a vector field with this kind of properties has on the global structure of a Lorentzian manifold,
when we have additional geometric information.

Let $(M,g)$ be a connected Lorentzian manifold with dimension $n>1$ and $U$ a closed reference frame. We will call $\Phi$ the flow of $U$,
$\omega$ the metrically equivalent one-form and $L_{p}$ the orthogonal leaf through $p\in M$. If there is not confusion we drop the point $p$
and write simply $L$.

\begin{proposition}
If $U$ is a complete and closed reference frame, then for all $p\in M$ the
map $\Phi:(\mathbb{R}\times L_p,-dt^{2}+g_{t})\rightarrow(M,g)$ is a normal
Lorentzian covering map, where $g_{t}$ is a metric tensor on $L_{p}$ for
each $t\in \mathbb{R}$ and $g_{0}=g_{\mid_{L_{p}}}$.
\end{proposition}

\begin{proof}
We have $L_{U}\omega =d\circ i_{U}\omega +i_{U}\circ d\omega =0$,
and therefore $\Phi $ is foliated, i.e. $\Phi _{t}(L_{p})=L_{\Phi
_{t}(p)}$ for all $t\in \mathbb{R}$ and $p\in M$. Thus, the map
$\Phi :\mathbb{R}\times L_{p}\rightarrow M$ is a local
diffeomorphism which is onto and $\Phi ^{\ast }g=-dt^{2}+g_{t}$,
where $g_t$ is a metric on $L_p$.

Now, we show that it is a covering map. Let $\sigma :[0,1]\rightarrow M$ be
a geodesic and $(t_{0},x_{0})\in \mathbb{R}\times L_{p}$ a point such that $%
\Phi (t_{0},x_{0})=\sigma (0)$. We must show that there exists a lift $%
\alpha :[0,1]\rightarrow \mathbb{R}\times L_{p}$ of $\sigma $ through $\Phi $
starting at $(t_{0},x_{0})$, \cite{ONeill}. There is a geodesic $\alpha
:[0,s_{0})\rightarrow \mathbb{R}\times L_{p}$, $\alpha (s)=(t(s),x(s))$,
such that $\Phi \circ \alpha =\sigma $ and $\alpha (0)=(t_{0},x_{0})$
because $\Phi $ is a local isometry. If we suppose $s_{0}<1$, there is a
geodesic $(t_{1}(s),x_{1}(s))$ such that $\Phi (t_{1}(s),x_{1}(s))=\sigma
(s) $ with $s\in (s_{0}-\varepsilon ,s_{0}+\varepsilon )$, then in the open
interval $(s_{0}-\varepsilon ,s_{0})$ it holds $\Phi (t(s),x(s))=\Phi
((t_{1}(s),x_{1}(s))$. Differentiating and using that $\Phi $ is foliated,
it is easy to see that $t_{1}(s)-t(s)=c\in \mathbb{R}$. Therefore it exists $%
\lim_{s\rightarrow s_{0}}\alpha (s)$ and the geodesic $\alpha $ is
extendible.

It remains to show that the group of deck transformations acts transitively
on the fibre. Take $(t_0,x_0)\in\mathbb{R}\times L_p$ such that $%
\Phi(t_0,x_0)=\Phi(0,p)=p$. Since $\Phi$ is foliated, it follows that $%
\Phi_{-t_0}(L_p)=L_p$ and thus the map $\mathbb{R}\times L_p\rightarrow
\mathbb{R}\times L_p$ given by $(t,x)\rightarrow (t+t_0,\Phi_{-t_0}(x))$ is
a deck transformation and takes $(0,p)$ to $(t_0,x_0)$.
\end{proof}

The following proposition, for the case that the flow is conformal when restricted on the leaves, is proved in a general form in \cite{Ponge},
but it does not include an explicit expression for the warping function. We include here a sketch of the proof to make the paper self-contained.

\begin{proposition}
\label{localdecomposition} Let $I\subset \mathbb{R}$ be an open interval, $L$
a manifold and $g$ a Lorentzian metric on $I\times L$ such that the
canonical foliations are orthogonal. If $U=\frac{\partial }{\partial t}$ is
an orthogonally conformal reference frame such that $\nabla div\,U$ is
proportional to $U$, then $g$ is the warped product $-dt^{2}+f(t)^{2}g_{0}$
where $g_{0}=g_{\mid _{L}}$ and $f(t)=\exp (\int_{0}^{t}\frac{div\,U(s,q)}{%
n-1}ds)$, being $q$ a fixed point of $L$.
\end{proposition}

\begin{proof}
 It is clear that $U$ is also orthogonally integrable, that is, the distribution $U^\perp$ is integrable, thus
$L_{U}g(v,w)=2\frac{div\,U}{n-1}g(v,w)$ for all $v,w\in U^{\perp }$. Therefore $\Phi _{t}:L_{(0,p)}\rightarrow L_{\Phi _{t}(0,p)}$ is a
conformal diffeomorphism with factor $\exp \left(
2\int_{0}^{t}\frac{div\,U(\Phi _{s}(0,p))}{n-1}ds\right) $. Since $%
L_{(t,p)}=\{t\}\times L$ and $\Phi _{t}(s,p)=(t+s,p)$, we have
\begin{equation*}
g((0_{t},v_{p}),(0_{t},w_{p}))=\exp \left( 2\int_{0}^{t}\frac{div\,U(s,p)}{%
n-1}ds\right) g_{0}(v,w)
\end{equation*}%
for all $v,w\in T_{p}L$. But $div\,U$ is constant on $L$ and the proposition follows.
\end{proof}

If a reference frame $U$ is closed, orthogonally conformal and $\nabla div\, U$ is proportional to $U$, we say that $U$ is a \textit{warped
reference frame}. This name is justified by the following corollary, which is a combination of the above two results.

\begin{corolario}
\label{covGRW} Let $(M,g)$ be a Lorentzian manifold and $U$ a complete
warped reference frame. Then the map $\Phi :(\mathbb{R}\times
L_{p},-dt^{2}+f(t)^{2}g_{0})\rightarrow (M,g)$ is a normal Lorentzian
covering for all $p\in M$, where $f(t)=\exp \left( \int_{0}^{t}\frac{%
div\,U(\Phi _{s}(p))}{n-1}ds\right)$ and $g_0=g|_{L_p}$.
\end{corolario}

In a warped product $-dt^2+f(t)^2g_0$ the vector field $\frac{\partial}{\partial t}$ is a warped reference frame. Observe that if we drop the
completeness hypothesis in the above corollary we do not obtain a covering map, but we do obtain a local
splitting around any point as $-dt^2+f(t)^{2}g_{0}$ with $U=\frac{%
\partial}{\partial t}$. Thus, the existence of a warped reference frame is
equivalent to the local splitting as a warped product $-dt^2+f(t)^2g_0$.

\begin{remark}\label{obs} \emph{If $M$ admits a complete warped reference frame then it is isometric to a quotient $\left( \mathbb{R}\times
_{f}L_{p}\right) /\Gamma $, where $\Gamma $ is a group of
isometries which preserves the canonical foliations. If $\psi \in
\Gamma $ then it is of the form $\psi
(t,x)=(t+k_{\psi },\varphi (x))$, for certain constant $k_{\psi }$, where $%
\varphi $
is an homothety of coefficient $c^{2}$ and $f(t)=cf(t+k_{\psi })$. Therefore, in order to show that $M$ is globally a GRW space we
must assure that $k_{\psi }=0$ for all $\psi \in \Gamma $, or equivalently, integral curves of $U$ meet the orthogonal leaves at only one value
of its parameter.}

\end{remark}

Recall that the quotient $\left( \mathbb{R}\times _{f}L_{p}\right) /\Gamma $
is not necessarily isometric to a warped product $\mathbb{S}^{1}\times _{f}N$
(called twisted product in \cite{Kuhnel}) even if $f$ is periodic, as the
following example shows.

\begin{example}
\label{ejemplo1}\emph{ Take $\tilde{M}=\mathbb{R}\times\mathbb{R}$ with the metric $g=-dt^2+(2+\cos(2\pi t))^2 dx^2$ and $\Gamma$ the group of
isometries generated by $\psi(t,x)=(t+1,x+1)$. Consider $M=\tilde{M}/\Gamma$. It has a complete warped reference frame $U$ which lifts to
$\frac{\partial}{\partial t}$, since $\frac{\partial}{\partial t}$ is invariant under $\Gamma$, but $M$ is neither a global GRW nor a warped
product of type $\mathbb{S}^1\times L$, see remark \ref{obs}.}
\end{example}

On the other hand, we can construct a quotient of a GRW that is not a global
GRW even if $f$ is non periodic.

\begin{example}
\label{ejemplo2}\emph{ Consider $\tilde{M}=\mathbb{R}\times \mathbb{R}^{2}$, with the metric $g=-dt^{2}+e^{t}(dx^{2}+dy^{2})$ which is a portion
of de Sitter spacetime $\mathbb{S}_1^3$, \cite{Alias2}. Let $\Gamma $ be the isometry group generated by $\psi
(t,(x,y))=(t+1,e^{-\frac{1}{2}}(x,y))$. It acts proper and discontinuously and the quotient $\tilde{M}/\Gamma$ is neither a global GRW nor a
warped product of type $\mathbb{S}^1\times_f L$.}
\end{example}

A timelike, closed and conformal vector field $V$ is characterized by the equation $\nabla_{X} V=a\ X$ for all vector field $X$, where $a$ is
certain function. Call $\lambda=|V|$ and $U=\frac{V}{\lambda}$. It is easy to show that $a=U(\lambda)$ and $\lambda$ is constant through the
orthogonal leaves. Moreover, $U$ is a warped reference frame. We write some useful formulas and facts.

\begin{lemma}

\label{conformalvector}Let $V$ be a timelike, closed and conformal vector
field, $\lambda=|V|$, $U=\frac{V}{\lambda}$ and $L$ an orthogonal leaf. Then

\begin{enumerate}

\item If $h$ is a constant function through the orthogonal leaves of $V$,
then $U(h)$ is also constant through the orthogonal leaves.

\item $\triangle\lambda=-U(U(\lambda))-(n-1)\frac{U(\lambda)^{2}}{\lambda}$.

\item $div\, U=(n-1)U(\ln \lambda)$.

\item If $v\in TL$ is an unit vector, then
\begin{eqnarray*}
Ric(v+U) &=&Ric_{L}(v)+\frac{n-2}{\lambda }\left( \frac{U(\lambda )^{2}}{%
\lambda }-U(U(\lambda ))\right) , \\
Ric(v) &=&Ric_{L}(v)+\frac{1}{\lambda }\left( U(U(\lambda ))+(n-2)\frac{%
U(\lambda )^{2}}{\lambda }\right)
\end{eqnarray*}%
\end{enumerate}
where $Ric_{L}$ is the Ricci tensor of $L$.
\end{lemma}

The above formulas can be directly proved from the fact that $V$ is closed and conformal, but observe that if we take $U$ the unitary of a
closed and conformal vector field $V$ in corollary \ref{covGRW}, then we get that $M$ is locally isometric to $\mathbb{R}\times _{f}L$, where
$f(t)=\frac{\lambda (\Phi _{t}(p))}{\lambda (p)}$ and we can deduce the above lemma from the standard formulas for a warped product,
\cite{ONeill}.

We finish this section establishing the relation between warped reference
frames and closed and conformal vector fields.

\begin{lemma}

\label{UnitarioConforme} Let $M$ be a complete Lorentzian manifold and $U$ a
warped reference frame. Then there exists a closed and conformal vector
field $V$ such that $U=\frac{V}{|V| }$.
\end{lemma}

\begin{proof}
\ Take $\Phi :\mathbb{R}\times _{f}L\rightarrow M$ the normal covering given in corollary \ref{covGRW}. If $\psi (t,x)=(t+k,B(x))$ is a deck
transformation, $\frac{f(t+k)}{f(t)}=c$ for all $t\in \mathbb{R}$, being $c$
certain constant, then using that $\mathbb{R}\times _{f}L$ is complete, \cite%
{Sanchez}
\begin{eqnarray*}
\infty &=&\int_{0}^{\infty }f(t)dt=\sum_{n=0}^{\infty
}c^{n}\int_{0}^{k}f(t)dt, \\
\infty &=&\int_{-\infty }^{0}f(t)dt=\sum_{n=0}^{\infty }c^{-n-1}\int_{0}^{k}f(t)dt.
\end{eqnarray*}

Thus $c=1$ and $f(t+k)=f(t)$ for all $t\in\mathbb{R}$. The vector field $f%
\frac{\partial}{\partial t}$ is closed and conformal and it is preserved by
the deck transformations. Thus there exists a closed and conformal vector
field $V$ on $M$ such that $U=\frac{V}{|V|}$.
\end{proof}

\section{Global GRW decompositions}

In \cite{Kuhnel} it was proved that a complete Lorentzian manifold with non negative constant scalar curvature which admits a non parallel
warped reference frame is isometric to a global GRW space. Now we get decomposition theorems on Lorentzian manifolds using Ricci curvature
hypothesis. Before, we need the following lemma.

\begin{lemma}

\label{TeoCompacta}Let $(M,g)$ be a non compact Lorentzian
manifold and $U$ a closed reference frame with compact orthogonal
leaves. Then $M$ is isometric to $(I\times L,-dt^{2}+g_{t})$,
where $I\subset\mathbb{R}$ is an open interval and $L$ a compact
Riemannian manifold.
\end{lemma}

\begin{proof}
\ Let $A$ be the domain of $\Phi$, $L$ an orthogonal leaf and $(a,b)\subset%
\mathbb{R}$ the maximal interval of $\mathbb{R}$ such that $(a,b)\times
L\subset A$. We claim that it is the maximal definition interval of each
integral curve with initial value on $L$. Indeed, suppose that $\Phi
_{t}(p_{0})$ is defined in $(a,b+\delta)$ for some $p_{0}\in L$. There is a $%
\eta \in \mathbb{R}\ $such that $(-\eta ,\eta )\times L_{\Phi
_{b}(p_{0})}\subset A$. Since $L_{\Phi _{b}(p_{0})}$ is compact, $\Phi _{-%
\frac{\eta }{2}}:L_{\Phi _{b}(p_{0})}\rightarrow L_{\Phi _{b-\frac{\eta }{2}%
}(p_{0})}$ is onto, and therefore a diffeomorphism. Now, for an arbitrary $%
p\in L$, $\Phi _{t}(p)$ can be defined in $(a,b+\eta )$ and we obtain a
contradiction.

If there were $(t_{0},p)\in(a,b)\times L$ such that $\Phi(t_{0},p)\in L$, $%
t_{0}\neq0$, then $M=\cup_{t\in\lbrack0,t_{0}]}\Phi_{t}(L),$ and it would be
compact. Therefore $\Phi:(a,b)\times L\rightarrow M$ is an injective map and
we obtain the desired result.
\end{proof}

\begin{theorem}\label{TeoRic2}\label{TeoRic3}

Let $M$ be a complete and non compact Lorentzian manifold with $n\geq 3$ and $U$ a non parallel warped reference frame. If one of the following
conditions is true
\begin{enumerate}
\item $Ric(U)\leq 0$,

\item $Ric(v)\geq0 $ for all $v\in U^{\perp}$,

\item $Ric(w)\geq 0$ for all lightlike vector $w$,
\end{enumerate}
then $M$ is globally a GRW space.
\end{theorem}
\begin{proof}
\ Take $V$ the closed and conformal vector field provided in lemma \ref{UnitarioConforme}. Then $\lambda=|V|$ is a non constant function because
$U$ is non parallel. Let $\gamma(t)$ be an integral curve of $U$.

$(1)$ If $\gamma$ returned to $L$ then $\lambda(\gamma(t))$ would be periodic, since $\Phi_t$ is foliated and $\lambda$ is constant through $L$.
Using $Ric(U)=-(n-1)\frac{UU(\lambda)}{\lambda}$, the hypothesis implies that $\frac{d\lambda(\gamma(t))}{\partial t}\geq 0$. Contradiction. The
result follows applying corollary \ref{covGRW}.

$(2)$ First suppose $\triangle\lambda\leq0$. Take $L$ an orthogonal leaf and suppose the $\gamma(0)\in L$. We know that if $\gamma$ returned to
$L$ then
$\lambda(\gamma(t))$ would be periodic. We have $%
-\triangle\lambda=U(U(\lambda ))+(n-1)\frac{U(\lambda)^{2}}{\lambda}$, thus if $z(t)=\ln\frac{ \lambda(\gamma(t))}{\lambda(\gamma(0))}$ we get
$0\leq z^{\prime\prime }+nz^{\prime2}$. Since $\lambda(\gamma(t))$ is periodic
and non constant there exists $t_{0}<t_{1}$ such that $z^{\prime}(t_{0})=z^{%
\prime}(t_{1})=0$ and $z^{\prime}(t)\neq0$ for all $t\in(t_{0},t_{1})$. Suppose that $z^{\prime }>0$ in $(t_{0},t_{1})$ (the case $z^{\prime}<0$
is similar). Then for all $t\in(t_{0}+\varepsilon,t_{1})$, where $\varepsilon$ is small enough, we get
\begin{equation*}
\int_{t_{0}+\varepsilon}^{t}\frac{-z^{\prime\prime}}{z^{\prime2}}\leq \int_{t_{0}+\varepsilon}^{t}n.
\end{equation*}

Thus $\frac{1}{z^{\prime}(t)}\leq n(t-t_{0}-\varepsilon )+\frac{1}{
z^{\prime}(t_{0}+\varepsilon)}$ and we get a contradiction taking $%
t\rightarrow t_{1}$. Therefore, the covering map $\Phi:\mathbb{R}\times
L\rightarrow M$ is injective (see remark \ref{obs}) and corollary \ref%
{covGRW} finishes the proof.

Now, suppose that there exists a point $p\in M$ with $\triangle\lambda(p)>0$ and call $L$ the leaf through $p.$ Take $z\in L$ and $v\in T_{z}L$
an unit vector. Since
\begin{equation*}
Ric_{L}(v)=Ric(v)-\frac{1}{\lambda}\left( U(U(\lambda))+(n-2)\frac {%
U(\lambda)^{2}}{\lambda}\right),
\end{equation*}
we get $Ric_{L}(v)\geq-\frac{1}{\lambda}\left( U(U(\lambda))+(n-1)\frac {
U(\lambda)^{2}}{\lambda}\right) =\frac{\triangle\lambda}{\lambda}(z).$ But $%
\lambda$ and $\triangle\lambda$ are constant through the orthogonal leaves,
thus $Ric_{L}(v)\geq\frac{\triangle\lambda}{\lambda}(p)>0$ for all unit vector $%
v\in TL.$ Using the completeness of $M$ and the theorem of Myers we conclude that $L$ is compact. Lemma \ref{TeoCompacta} and proposition
\ref{localdecomposition} show that $M$ is globally a GRW space.

$(3)$ Let $h:M\rightarrow\mathbb{R}$ be given by $h=U(U(\lambda ))-\frac{U(\lambda)^{2}}{\lambda}$. We consider two possibilities: there exists
$p\in M$ with $h(p)>0$ or $h(q)\leq0$ for all $q\in M$.

Assume the second one. The function $z(t)=\ln\frac{\lambda(\gamma(t))}{ \lambda(\gamma(0))}$ verifies $z^{\prime\prime}\leq0$. Therefore
$\gamma$ can not return to $L_{\gamma(0)}$ and applying corollary \ref{covGRW}, $M$ is globally a GRW space.

Suppose now that there exists $p\in M$ such that $h(p)>0$ and take $L$ the leaf through $p.$ If $v\in T_{z}L$ is an unit vector then $v+U_{z}$
is lightlike and
\begin{equation*}
0\leq Ric(v+U_{z})=Ric_{L}(v)+\frac{n-2}{\lambda}\left( \frac{ U(\lambda)^{2}%
}{\lambda}-U(U(\lambda))\right).
\end{equation*}

Then we get $\frac{n-2}{\lambda(z)}h(z)\leq Ric_{L}(v).$ But $\lambda$ and $%
h $ are constant through the orthogonal leaves. Hence $0<\frac{n-2}{\lambda (p)}h(p)\leq Ric_{L}(v)$ for all unit vector $v\in TL$ and we can
conclude as in the point 2.
\end{proof}

We can not use $Ric(U)\geq 0$ because in this case $U$ would be parallel. If $Ric(u)\geq0$ for all timelike vector $u$ then it is said that the
timelike convergence condition holds (TCC) and if $Ric(u)\geq0$ for all lightlike vector $u$ then it is said that the null convergence condition
holds (NCC). We can not suppose the more restrictive TCC condition because $U$ would be parallel too, \cite{Sanchez}.

A condition like $div\,U\geq 0$ (resp. $div\,U\leq 0$), leads trivially to a GRW decomposition because $\lambda(\gamma(t))$ would be increasing
(resp. decreasing).

\begin{corollary}

Let $(L,g_{0})$ be a non compact and complete Riemannian manifold and $M=%
\mathbb{S}^{1}\times L$ endowed with a warped product metric $%
-dt^{2}+f(t)^{2}g_{0}$. If the NCC holds then $f$ is constant.
\end{corollary}

Causality hypotheses are frequently used in Lorentzian geometry, besides curvature hypotheses. Since the injectivity of $\Phi$ depends on the
behavior of some timelike curves it seems natural to impose a causality condition to reach global decompositions. However, a hard condition like
being globally hyperbolic is not sufficient to obtain a global product, as the following example shows.

\begin{example}
\emph{(Compare with proposition 2 in \cite{Harris2}) Take $\tilde{M}=\mathbb{R}%
^{2} $ with the metric $-dt^{2}+f(t)^{2}dx^{2}$, where $f(t)=4+\sin(2\pi t)$%
. Call $\Gamma$ the isometry group generated by $\Psi (t,x)=(t+1,x+1)$ and $%
\Pi:\tilde{M}\rightarrow M=\tilde{M}/\Gamma$ the projection. The vector field $V=\Pi_{*}(f\frac{\partial}{\partial t})$ is timelike, closed and
conformal. The manifold $M$ verifies any causality condition. In fact, $\Pi \left(\{(t,x):t=x\}\right)$ is a Cauchy hypersurface. But $M$ does
not split as a GRW, see remark \ref{obs}.}
\end{example}

\begin{lemma}

Let $M$ be a complete Lorentzian manifold and $V$ a timelike, closed and
conformal vector field with unitary $U$. If there exists $t_0>0$ and $p\in M$
such that $\Phi_{t_0}(p)\in L_p$, being $\Phi$ the flow of $U$, then

\begin{enumerate}

\item $\Phi_{t_0}:M\rightarrow M$ is an isometry.

\item If $M$ is chronological, then the isometry group $\Omega$
generated by $\Phi_{t_0}$ is isomorphic to $\mathbb{Z}$ and acts
on $M$ in a properly discontinuous manner.
\end{enumerate}

\end{lemma}

\begin{proof}
\ $(1)$ Use that $V$ is conformal, $\Phi_t$ is foliated and $\lambda=|V|$ is
constant through the orthogonal leaves.

$(2)$ Suppose that the manifold is future oriented by the vector field $U$.
We will construct for each $q\in M$ an open set $\Theta$ of $q$ such that $%
\Phi_{nt_0}(\Theta)\cap \Theta=\emptyset$ for all $n\in \mathbb{Z}-\{0\}$. Recall that $\Phi:\mathbb{R}\times _f L_q \rightarrow M$ is a local
isometry with $f(t)=\frac{\lambda(\Phi_q(t))}{\lambda(q)}$, see comments after lemma \ref{conformalvector}. Take $k=max\{f(t)^2\, :
\, t\in [-t_0,t_0]\}$ and the ball $W=B(q,\frac{\varepsilon}{2\sqrt{k}}%
)\subset L_q$, where $\varepsilon<t_0$ is small enough to $%
\Phi:(-\varepsilon,\varepsilon)\times_f W\rightarrow \Theta$ be an isometry
and $W$ a normal neighborhood of $q$. Given $\Phi_s(z)\in \Theta$, with $%
z=exp_q(v)\in W$ and $s\in (-\varepsilon,\varepsilon)$, we can construct
future timelike curves for $t \in [0,1]$
\begin{eqnarray*}
\alpha(t)&=&\Phi\left(\frac{\varepsilon}{2}(t-1),\exp_q{((1-t)v)}\right) \\
\beta(t)&=&\Phi\left(\frac{\varepsilon}{2}t,\exp_q{(tv)}\right)
\end{eqnarray*}
from $\Phi_{-\frac{\varepsilon}{2}}(z)$ to $q$ and from $q$ to $\Phi_{\frac{%
\varepsilon}{2}}(z)$ respectively.

Suppose $\Phi _{nt_{0}}(\Theta )\cap \Theta \neq \emptyset $, $n\in \mathbb{Z%
}-\{0\}$. Take $x,y\in \Theta $ with $y=\Phi _{nt_{0}}(x)$. The open set $%
\Theta $ is an union of segments $\Phi _{z}(-\varepsilon ,\varepsilon )$ of
integral curves of $U$, with $z\in W$, thus we have $x\in \Phi
_{z_{1}}(-\varepsilon ,\varepsilon )$, $y\in \Phi _{z_{2}}(-\varepsilon
,\varepsilon )$. We can use a timelike curve $\alpha $ from $\Phi
_{-(\varepsilon /2)}(z_{2})$ to $q$ and $\beta $ from $q$ to $\Phi
_{(\varepsilon /2)}(z_{1})$. The curve formed by the segment of  $\Phi
_{z_{1}}$ from $\Phi _{(\varepsilon /2)}(z_{1})$ to $\Phi _{-(\varepsilon
/2)}(z_{2})$, $\alpha $ and $\beta $, is timelike and closed in
contradiction with the chronology hypothesis. Thus, $\Phi _{nt_{0}}(\Theta
)\cap \Theta =\emptyset $ for all $n\in \mathbb{Z}-\{0\}$.

Take the Riemannian metric $g_{R}=g+2\omega \otimes \omega$. The
group $\Omega$ is a group of isometries for $g_{R}$ too. Since
$g_{R}$ is Riemannian, the existence of the above open set
$\Theta$ for each $q\in M$ is sufficient to show that the action
of $\Omega$ in $M$ is properly discontinuous.
\end{proof}

Note that we must use a Riemannian argument in the above proof because the analogous statement in Lorentzian geometry is not true.

\begin{theorem}

Let $M$ be a chronological complete Lorentzian manifold and $U$ a non parallel warped reference frame. Then $M$ is a global GRW space or there
is a
Lorentzian covering map $\Psi :M\rightarrow \mathbb{S}^{1}\times N$, where $(%
\mathbb{S}^{1}\times N,-dt^{2}+f^{2}g_{N})$ is a Lorentzian warped product.
\end{theorem}

\begin{proof}
\ Suppose $M$ is not a GRW space. Then there is $p\in M$ such that the
integral curve $\Phi_{p}(t)$ of $U$ returns to $L$, the leaf of $U^{\perp}$
through $p$ (remark \ref{obs}). Take $t_{0}=\inf \{t>0:\Phi_{p}(t)\in L\}$
and call $\lambda(t)=\lambda(\Phi_{p}(t))$, where $\lambda=|V|$ and $V$ is
the closed and conformal vector field with $U=\frac{V}{|V|}$ (lemma \ref%
{UnitarioConforme}). We can suppose that $\lambda^{\prime}(0)>0$, since $%
\lambda(t)$ is not constant. Then, $t_{0}>0$ and it is a minimum. Therefore $%
\Phi_{p}(t_{0})\in L$ and applying the above lemma, $\Omega$ acts on $M$ in a properly discontinuous manner. Using corollary \ref{covGRW} it is
easy to show that $M/\Omega$ is isometric to a Lorentzian warped product $(\mathbb{S}^{1}\times N,-dt^2+f^{2}g_{N})$, where $N=L/\Omega$.
\end{proof}

\section{Perfect Fluids}

A Lorentzian four dimensional manifold $M$ is called a perfect fluid if there is a reference frame $U$ and $\rho ,\eta \in C^{\infty }(M)$ (the
energy and pressure) such that the stress-energy tensor is $T=(\rho +\eta )\omega \otimes \omega +\eta g$, or equivalently $Ric=(\rho +\eta
)\omega\otimes \omega+\frac{1}{2}(\rho -\eta )g$. A perfect fluid satisfies $U(\rho)=-(\rho +\eta)div\,U$ (energy equation ) and $(\rho
+\eta)\nabla_UU=-\nabla^{\perp}\eta$ (force equation), where $-\nabla^{\perp}\eta$ is the component of $\nabla\eta$ orthogonal to $U$.

It is well known that any Robertson-Walker (RW) spacetime is a perfect fluid \cite{ONeill}, and it is straightforward to see that a four
dimensional GRW spacetime is a perfect fluid if and only if it is a RW spacetime.

Moreover, in a RW spacetime we have the following basic relation between the warping function $f$ and pressure and energy.
\begin{equation*}
\frac{3f''}{f}=-\frac{1}{2}(\rho +3\eta),\ \  \rho'=-3(\rho +\eta)\frac{f'}{f}.
\end{equation*}

Eliminating $f$, the above equations implies the following equation of state in the open set $\rho +\eta \neq 0$.

\begin{equation}\label{FP}
\left( \frac{\rho^{\prime}}{\rho+\eta }\right) ^{\prime}=\frac{1}{3}\left( \frac{\rho^{\prime}}{\rho+\eta}\right) ^{2}+\frac{1}{2}(\rho+3\eta).
\end{equation}

A perfect fluid is called barotropic if it satisfies an equation of state $\eta =\eta (\rho )$. Observe that a RW perfect fluid
$\mathbb{R}\times _{f}L$ is barotropic in the open set $\rho +\eta \neq 0$ and $\frac{d\rho }{d\eta}\neq 0$.

A natural question is which conditions on a perfect fluid make it a global Robertson-Walker space, \cite{Garcia3,Garcia2}. Using the
decomposition theorems of the previous section we can give an answer to this question. In fact, the following theorem shows that under mild
conditions, the equation of state (\ref{FP}) give rise to a global decomposition as a RW space.

\begin{theorem}

Let $M$ be a (four dimensional) non compact spacetime with a barotropic perfect fluid $(U,\rho,\eta)$ such that $U$ is geodesic,
$\frac{d\eta}{d\rho}\neq0$, $\frac{ d^{2}\eta }{d^{2}\rho}\neq0$, $\rho+\eta\neq 0$ and $\rho>0$ is not constant. If it verifies the equation of
state (\ref{FP}), then either $M$ is incomplete or a global Robertson-Walker space.
\end{theorem}

\begin{proof}
\ We first show that $U$ is closed, that is, $d\omega =0$. From the force equation $g(X,\nabla \eta )=-(\rho +\eta )g(X,\nabla _{U}U)$ for all
$X\in U^{\perp }$ we get $d\eta =h\,\omega $, being $h$ certain function, and
\begin{equation*}
-X(h)=X(g(U,\nabla \eta ))=g(U,\nabla _{X}\nabla \eta )=g(X,\nabla
_{U}\nabla \eta )=0.
\end{equation*}

Hence $0=dh\wedge \omega+h\,d\omega=h\,d\omega$. Thus, to show $%
d\omega_{p}=0 $ it is enough to consider a critical point $p$ of $\eta$.
Call $\eta(t)=\eta(\gamma(t))$ and $\rho(t)=\rho(\gamma(t))$ where $\gamma$
is the integral curve of $U$ with $\gamma(0)=p$. Then $\eta^{\prime}=0$ if
and only if $\rho^{\prime }=0$. If $\rho^{\prime}(t)=0$ for all $%
t\in(-\varepsilon,\varepsilon)$ then $\rho+3\eta=0$ in contradiction with $%
\frac{d^{2}\eta}{d^{2}\rho}\neq0$. Thus there is a sequence $t_{n}$
converging to $0$ with $\rho^{\prime}(t_{n})\neq0 $ and therefore $%
d\omega_{p}=0$.

Now, consider the self-adjoint endomorphism $A:U^{\perp }\rightarrow U^{\perp }$ given by $A(X)=\nabla _{X}U$. A straightforward computation
give us $Ric(U)=-U(div\,U)-\Vert A\Vert ^{2}$, (to prove this take $\{U_{1},U_{2},...,U_{n}\}$ a local orthonormal basis around an arbitrary
point $p$, such that $\nabla _{U_{i}}U_{j}(p)=0$ and $U_{1}(p)=U(p)$). From this, the energy equation and the equation of state (\ref{FP}), we
get
\begin{equation*}
\Vert A\Vert ^{2}=\left( \frac{\rho ^{\prime }}{\rho +\eta
}\right) ^{\prime }-\frac{1}{2}(\rho +3\eta )=\frac{1}{3}\left(
\frac{\rho ^{\prime }}{\rho +\eta }\right) ^{2}=\frac{1}{3}(trace\
A)^{2}.
\end{equation*}

Then $A=\frac{div\, U}{3}id$ and $U$ is orthogonally conformal. We know that $\eta$ and $\rho$ are constant on the orthogonal leaves. The energy
equation implies that $div\,U$ is constant through the orthogonal leaves, that is, $\nabla div\, U$ is proportional to $U$ and hence $U$ is a
warped reference frame. If $U$ were parallel then $M$ would be locally a direct product, which implies that $\rho$ is constant, \cite{ONeill}.

If $M$ is complete, we apply theorem \ref{TeoRic2}(2) if $\rho +\eta <0$ and theorem \ref{TeoRic3}(3) if $\rho +\eta >0$. (See the expression of
the Ricci tensor of a perfect fluid at the beginning of this section).
\end{proof}

\section{Lightlike sectional curvature and timelike hypersurfaces}

Let $M$ be a Lorentzian manifold and $U$ a reference frame. We can
define a curvature for degenerate planes as follows,
\cite{Harris1}. Take $\Pi$ a degenerate plane and a basis
$\{u,v\}$, where $u$ is the unique lightlike vector in $\Pi$ with
$g(u,U)=1$. The lightlike sectional curvature of $\Pi$ is
\begin{equation*}
\mathcal{K}_{U}(\Pi)=\frac{g(R(v,u,u),v)}{g(v,v)}.
\end{equation*}

This sectional curvature depends on the choice of the reference frame $U$,
but its sign does not change if we choose another vector field. In fact, if $%
U^{\prime}$ is another reference frame then $\mathcal{K}_{U}(\Pi)=g(u,U^{%
\prime})^{2} \mathcal{K}_{U^{\prime}}(\Pi)$ where $u\in\Pi$ is the unique
lightlike vector such that $g(u,U)=1$. Thus, it makes sense to say positive
lightlike sectional curvature or negative lightlike sectional curvature.

If $U$ is geodesic, $n>3$ and $\mathcal{K}_{U}$ is a never zero point function then $U$ is a warped reference frame \cite{Harris2,Koch} and,
under completeness hypothesis, it follows from corollary \ref{covGRW} that $M$ is covered by a Robertson-Walker space $\mathbb{R}\times L$. As
an application of the decomposition results presented above, we obtain conditions on the lightlike sectional curvature which ensure the global
decomposition of $M$. First, we give some relations between lightlike sectional curvature and Ricci curvature.

\begin{lemma}

\label{Ricciycurvnula} Let $M$ be a Lorentzian manifold and $U$ a timelike and unit vector field. If $u$ is a lightlike vector with $g(u,U)=1$,
then
\begin{equation*}
Ric(u)=\sum_{i=1}^{n-2}\mathcal{K}_{U}(span(u,e_{i})),
\end{equation*}
where $\{e_{1},...,e_{n-2}\}$ is an orthonormal basis of $u^{\perp}\cap
U^{\perp}$.
\end{lemma}

\begin{lemma}

\label{formulacurvnula} Let $M=I\times _{f}L$ be a GRW space and $u$ a
lightlike vector such that $g(u,\frac{\partial }{\partial t})=1$. If $u=-%
\frac{\partial }{\partial t}+w$, $w\in TL$, and $\Pi =span(u,v)$ is a
degenerate plane, where $v\in TL$, then $\mathcal{K}_{\frac{\partial }{%
\partial t}}(\Pi )=\frac{K^{L}(span(v,w))+f^{\prime 2}-f^{\prime \prime }f}{%
f^{2}}$, where $K^{L}$ is the sectional curvature of $L$. If moreover $L$
has constant sectional curvature then $Ric(u)=(n-2)\mathcal{K}_{\frac{%
\partial }{\partial t}}(\Pi )$.
\end{lemma}

\begin{proof}
\ It is a straightforward calculation. If $L$ has constant sectional
curvature, then $\mathcal{K}_{\frac{\partial}{\partial t}}$ is a point
function and the above lemma gives us $Ric(u)=(n-2)\mathcal{K}_{\frac{%
\partial}{\partial t}}(\Pi)$.
\end{proof}

\begin{lemma}

\label{lemaincompl} Let $M=\mathbb{R}\times_{f}L$ be a GRW space. If $M$ is
lightlike complete and $Ric(u)>0$ for all lightlike vector $u$, then $%
Ric_{L}(v)> 0$ for all $v\in TL$.
\end{lemma}

\begin{proof}
\ Suppose there is $v\in TL$ unit vector such that $Ric_{L}(v)\leq 0$. Then
\begin{equation*}
0<Ric(\frac{\partial}{\partial t}+v)=Ric_{L}(v)+\frac{n-2}{f}\left( \frac{%
f^{\prime2}}{f}-f^{\prime\prime }\right).
\end{equation*}

If we call $y=\ln f$ then $y^{\prime\prime}=\frac{f^{\prime\prime}f-f^{%
\prime 2}}{f^{2}}<0$. We can suppose $y^{\prime}(0)>0$. Now, $%
\int_{-\infty}^{0}e^{y}\leq\int_{-\infty}^{0}e^{y^{\prime}(0)t+y(0)}<\infty$
and we conclude from \cite{Sanchez} that $M$ is lightlike incomplete.
\end{proof}

\begin{proposition}

\label{TeoKpositiva} Let $M$ be a non compact and complete Lorentzian manifold with $n>3$ and $U$ a geodesic reference frame. If the lightlike
sectional curvature $\mathcal{K}_{U}$ is a positive point function, then $M$ is globally a Robertson-Walker space $\mathbb{R}\times L$ with $L$
of constant positive sectional curvature.
\end{proposition}

\begin{proof}
\ As we said at the beginning of this section, $M$ is covered by a RW space $%
\mathbb{R}\times L$. From lemmas \ref{Ricciycurvnula}, \ref{formulacurvnula} and \ref{lemaincompl}, $L$ has positive constant curvature and
therefore is compact. Now, we apply lemma \ref{TeoCompacta}.
\end{proof}

\begin{theorem}

Let $M$ be a non compact and complete Lorentzian manifold with $n>3$ and $U$
a geodesic reference frame. Suppose that the lightlike sectional curvature $%
\mathcal{K}_{U}$ is a never zero point function such that $\frac{1}{n-2}%
Ric(U)<\mathcal{K}_{U}$. Then $M$ is globally a Robertson-Walker
space.
\end{theorem}

\begin{proof}
\ If $U$ were parallel then $Ric(U)=0$ and the above proposition give us the desired result. If $U$ is not parallel, take $v,w\in U^{\perp}$
unit and orthogonal vectors and the degenerate plane $\Pi=span(-U+w,v)$. Since $M$ is
locally a RW space, we can apply lemma \ref{formulacurvnula} and we get $%
Ric(-U+w)=(n-2)\mathcal{K}_{U}(\Pi)$. But $Ric(-U+w)=Ric(U)+Ric(w)$ and
therefore $Ric(w)>0$ for all $w\in U^{\perp}$. Now, we apply theorem \ref%
{TeoRic2}(2).
\end{proof}

Spacelike hypersurfaces are widely studied in General Relativity due to
their role as initial data hypersurfaces in the Cauchy problem. On the other
hand, a k-dimensional timelike submanifold can be interpreted as the history
of a (k-1)-dimensional spacelike submanifold. Timelike totally umbilic
hypersurfaces are called photon surfaces and were studied in \cite{Ellis}.
Now, we show that in a local GRW spaces, photon surfaces are global GRW
spaces under certain curvature condition.

\begin{lemma}

Let $M$ be a Lorentzian manifold and $U$ a reference frame. If $S$ is a
timelike totally umbilic hypersurface and $u\in TS$ is a lightlike vector
with $g(u,U)=1$, then
\begin{equation*}
Ric_{S}(u)=\sum_{i=1}^{n-3}\mathcal{K}_{U}(span(u,e_{i})),
\end{equation*}
where $\{e_{1},...,e_{n-3}\}$ is an orthonormal basis of $u^{\perp}\cap
U^{\perp}\cap TS$.
\end{lemma}

\begin{proof}
\ Suppose that $II(X,Y)=g(X,Y)N$ for all $X,Y\in TS$. Take $u\in TS$ a
lightlike vector such that $g(u,U)=1$ and $\{v_{1},...,v_{n}\}$ an
orthonormal basis with $v_{n}$ the unitary of $N$. By the Gauss equation and
the fact that $u$ is lightlike
\begin{equation*}
Ric(u)=\sum_{i=1}^{n} \varepsilon_{i}g(R(v_{i},u,u),v_{i})=Ric_{S}(u)+%
\mathcal{K}_{U}(span(u,v_{n})).
\end{equation*}

Now, take $\{e_{1},...,e_{n-2}\}$ an orthonormal basis of $u^{\perp}\cap
U^{\perp}$ such that $e_{n-2}$ is the unitary of the projection of $N$ on $%
u^{\perp}\cap U^{\perp}$. Then $span(u,v_{n})=span(u,e_{n-2})$ and using lemma \ref{Ricciycurvnula},
$Ric(u)=\sum_{i=1}^{n-2}\mathcal{K}_{U}(span(u,e_{i}))$. Therefore $Ric_{S}(u)=\sum_{i=1}^{n-3}\mathcal{K}_{U}(span(u,e_{i}))$.
\end{proof}

\begin{theorem}

\label{umbilicalhypersurface} Let $M$ be a complete Lorentzian manifold with
$n\geq 4$ and $U$ a warped reference frame. If $S$ is a timelike, non
compact, complete and totally umbilic hypersurface of $M$ with never zero
mean curvature such that $\mathcal{K}_U(\Pi )\geq 0$ for all degenerate plane $%
\Pi$ tangent to $S$, then $S$ is globally a GRW space.
\end{theorem}

\begin{proof}
\ Take $V$ closed and conformal with $U=\frac{V}{|V|}$, (lemma \ref{UnitarioConforme}). Suppose $\nabla V=a\cdot id$ and $II(X,Y)=g(X,Y)N$ for
all $X,Y\in TS$, where $N$ is the
mean curvature vector field of $S$. If $V=\alpha \frac{N}{|N|}+W$, where $%
W\in TS$, then it follows that $\nabla_{X}^{S}W=(a+g(N,V))X$ for all $X\in TS $, i.e. $W$ is closed and conformal in $S$.

Suppose first that $W$ is parallel in $S$, that is, $a=-g(N,V)$. Call $%
c=g(W,W)$ and take the function on $S$ given by $h(p)=g(V,V)$
which is
constant on the orthogonal leaves of $W$. On the other hand, $W(h)=2ac$ and $%
W(a)=-W(g(V,N))=cg(N,N)+a\frac{W(g(N,N))}{2g(N,N)}$. If $\gamma$ is an integral curve of $W$ and $h(t)=h(\gamma (t))$, then $h^{\prime \prime
}(t_{0})>0$ at each critical point $t_{0}$ of $h$, showing that $h$ can not
be periodic. Therefore remark \ref{obs} and corollary \ref{covGRW} gives us that $%
S $ is globally a direct product.

Suppose now that $W$ is not parallel. The above lemma shows that $S$ satisfies the NCC condition, and theorem \ref{TeoRic3}(3) finishes the
proof.
\end{proof}

\begin{corollary}

\label{coroumbilical} Let $M$ be a complete Lorentzian manifold with $n\geq
4 $ and $U$ a parallel reference frame. If $K(\Pi)\geq 0$ for all plane $%
\Pi\in U^\perp$, then any timelike, non compact, complete and
totally umbilic hypersurface of $M$ with never zero mean curvature
is globally a GRW space.
\end{corollary}

\begin{proof}
\ If $\Pi=span(-U+v,w)$ is a degenerate plane to the hypersurface,
from lemma \ref{formulacurvnula} we have
$\mathcal{K}_{U}(\Pi)=K(span(v,w))\geq 0$.
\end{proof}

In general, it is a difficult question to check the completeness
of a timelike hypersurface $S$ of a complete Lorentzian manifold
$M$. In the case of totally umbilic hypersurfaces of a GRW space
we can give a simple criterium.

If $U$ is a reference frame on $M$, we can define the hyperbolic angle $%
\theta\in [0,\infty)$ between $S$ and $U$ as the hyperbolic angle between $U$ and the projection of $U$ in $S$. More explicitly, if $N_{0}$ is
the normal unit vector field to $S$ and $U=\alpha N_{0}+W,\, W\in N_{0}^{\perp}$, then $\theta$ is characterized by
\begin{equation*}
\cosh \theta=\frac{-g(U,W)}{\sqrt{-g(W,W)}}.
\end{equation*}

\begin{proposition}

Let $\mathbb{R}\times_{f} L$ be a complete GRW space. If $S$ is a timelike, closed (as a subset of $\mathbb{R}\times L$) and totally umbilic
hypersurface of $\mathbb{R}\times_{f} L$ such that the hyperbolic angle between $\frac{\partial}{\partial t}$ and $S$ is bounded, then $S$ is
complete.
\end{proposition}

\begin{proof}
\ Let $N_{0}$ be the normal unit vector field to $S$ and $V=f\frac{\partial }{\partial t}$. Put $II(X,Y)=g(X,Y)N$ and $V=\alpha N_{0}+W,\,W\in
N_{0}^{\perp }$.
We already know that $W$ is closed and conformal in $S$ and thus $U=\frac{W}{%
|W|}$ is a warped reference frame. The Riemannian metric $%
g_{R}=g+2dt^{2}=dt^{2}+f^{2}g_{0}$ is complete, \cite{ONeill,Sanchez},
and we have $g_{R}(U,U)=-1+2g(\frac{\partial }{\partial t}%
,U)^{2}=-1+2\cosh^{2}\theta $, where $\theta $ is the hyperbolic angle between $\frac{\partial }{\partial t}$ and $S$. Therefore $g_{R}(U,U)$ is
bounded and thus it is a complete vector field. Using corollary \ref%
{covGRW} there is a Lorentzian covering $\mathbb{R}\times _{h}Q\rightarrow S$%
, where $h(s)=\frac{|W|_{\gamma (s)}}{|W|_{\gamma (0)}}$ and $\gamma $ is an
integral curve of $U$. We can suppose without loss of generality $%
|W|_{\gamma (0)}=1$. Since $S$ is totally umbilic, its lightlike geodesics
are geodesics of $\mathbb{R}\times L$, thus $S$ and $\mathbb{R}\times _{h}Q$
are lightlike complete and hence $Q$ is complete. If we show that
\begin{equation*}
\int_{0}^{\infty }\frac{h(s)}{\sqrt{1+h(s)^{2}}}ds=\int_{-\infty }^{0}\frac{%
h(s)}{\sqrt{1+h(s)^{2}}}ds=\infty ,
\end{equation*}%
then $\mathbb{R}\times _{h}Q$ is complete and so is $S$, \cite{Sanchez}. Let $T:\mathbb{R}\times L\rightarrow \mathbb{R}$ be the projection and
consider the diffeomorphism $\rho =T\circ \gamma :\mathbb{R}\rightarrow \mathbb{R}$. Then $\frac{d}{ds}\rho =-g(\frac{\partial }{\partial
t},U)=\cosh \theta \leq c$, where $c$ is a certain constant. Since $g(W,W)_{\gamma (s)}\leq g(V,V)_{\gamma (s)}=-f^{2}(T(\gamma (s))$ it follows
that
\begin{equation*}
\int_{0}^{\infty }\frac{f\circ \rho }{\sqrt{1+(f\circ \rho
)^{2}}}ds\leq \int_{0}^{\infty }\frac{h}{\sqrt{1+h^{2}}}ds,
\end{equation*}%
but
\begin{equation*}
\infty =\int_{0}^{\infty }\frac{f(t)}{\sqrt{1+f(t)^{2}}}dt\leq
c\int_{0}^{\infty }\frac{f\circ \rho (s)}{\sqrt{1+(f\circ \rho
(s))^{2}}}ds.
\end{equation*}

\end{proof}

\begin{example}

\emph{Take $M=\mathbb{S}^{1}\times\mathbb{R}\times\mathbb{S}^{2}$ with metric $%
g=-dt^{2}+f(t)^{2}g_{0}$ where $f(t)=2+\cos t$, $g_{0}=dx^{2}+\frac{1}{2}%
g_{S}$ and $g_{S}$ is the canonical metric in $\mathbb{S}^2$. Take the
universal covering $\varrho:\tilde{M}\rightarrow M$, $F:\tilde{M}\rightarrow%
\mathbb{R}$ given by $F(t,x,p)=x-\int_{0}^{t}\frac{1}{2f(s)}ds$ and $\tilde{S%
}=F^{-1}(0)$. It follows that $\nabla F=\frac{1}{2f}\frac{\partial}{\partial
t}+\frac{1}{f^{2}}\frac{\partial}{\partial x}$ and $(\nabla F)^{\perp}=T%
\mathbb{S}^{2}\oplus span( \frac{\partial}{\partial t}+\frac{1}{2f}\frac{%
\partial}{\partial x})$. Now, a direct computation gives us that $%
Hess_{F}(X,Y)=\frac{f^{\prime}}{2f^{2}}g(X,Y)$ for all $X,Y\in\nabla
F^{\perp}$. Hence $\tilde{S}$ is a timelike, closed, non compact and totally
umbilic hypersurface of $\tilde{M}$. Since the hyperbolic angle between $%
\tilde{S}$ and $\frac{\partial}{\partial t}$ is constant, it follows that $%
\tilde{S}$ is complete.}

\emph{Take now $\Pi$ a degenerate plane tangent to $\tilde{S}$. Suppose $%
\Pi=span(u,v)$ where $g(u,u)=g(u,v)=g(\frac{\partial}{\partial t},v)=0$ and $%
g(v,v)=g(u,\frac{\partial}{\partial t})=1$. Then $u=-\frac{\partial}{%
\partial t}+w$, where $g(w,w)=1$ and $g(w,\frac{\partial}{\partial t}%
)=g(v,w)=0$.}

\emph{Since $\Pi $ is tangent to $\tilde{S}$, the structure of $(\nabla F)^{\perp
} $ implies that $v\in T\mathbb{S}^{2}$ and the relation $g(u,\nabla F)=0$
that $w=-\frac{1}{2f}\frac{\partial }{\partial x}+X$, where $X\in T\mathbb{S}%
^{2}$. Now, if we denote $K^{0}$ the sectional curvature of $(\mathbb{R}%
\times \mathbb{S}^{2},g_{0})$, then $K^{0}(span(v,w))=\frac{3}{2}$, therefore lemma \ref{formulacurvnula} give us}
\begin{equation*}
\mathcal{K}_{\frac{\partial }{\partial t}}(\Pi )=\frac{K^{0}(span(v,w))+f^{%
\prime 2}-f^{\prime \prime }f}{f^{2}}=\frac{\frac{5}{2}+2\cos (t)}{f^{2}}>0.
\end{equation*}%
\emph{The hypersurface $S=\varrho (\tilde{S})$ is not compact. Applying theorem \ref{umbilicalhypersurface} it is a global GRW space.}
\end{example}

\begin{example}
\emph{\ From corollary \ref{coroumbilical} we can deduce the well known fact that the anti de Sitter space $\mathbb{H}^n_{1}(r)=\{p\in
\mathbb{R}^{n+1}_{2} \, : \, \langle p,p\rangle=-r^{2}\}$ can not be immersed completely in Minkowski space $\mathbb{R}^{n+1}_{1}$. In fact, if
$S\subset\mathbb{R}^{n+1}_{1}$ is a complete Lorentzian manifold with constant negative curvature, then $S$ is totally umbilic, \cite{ONeill}.
On the other hand, the Gauss equation implies that it has never zero mean curvature vector. Using the corollary, $S$ is globally a GRW space.
But a GRW space of constant negative curvature is incomplete, \cite{Sanchez}.}
\end{example}

\section{Uniqueness of GRW decompositions}

In \cite{Esche}, the uniqueness of direct product decompositions of a
Riemannian manifold is studied, where only products into indecomposable
factors are considered and the uniqueness is understood in the following
sense: a direct product decomposition into indecomposable factors is unique
if the corresponding foliations are uniquely determined. Euclidean space
admits more than one direct product decomposition, but this is essentially
the only Riemannian manifold with this property.

We know that in a GRW decomposition $\mathbb{R}\times_f L$  the vector field $f\frac{\partial}{\partial t}$ is timelike, closed and conformal.
Conversely, a timelike, closed and conformal vector field gives rise to a warped reference frame and therefore to a local GRW decomposition.
Thus, to deal with the uniqueness of GRW decompositions, it is sufficient to study how many timelike, closed and conformal vector fields can
exists on a Lorentzian manifold, up to scalar multiplication.

The uniqueness question has been recently analyzed in \cite{Sanchez2} for static spacetimes. A possible interest of this kind of results for GRW
spacetimes comes from the recently introduced "big rip" models which try to explain the accelerated expansion of the Universe. The qualitative
properties of the models depend on the behavior of the warping function $f$, \cite{Fernandez}. So we must ensure that the qualitative behavior
of this function (or the function itself) is univocally determined.

\begin{example}\label{seudoesfera}\emph{ The de Sitter space
 $\mathbb{S}_{1}^{n}(r)=\{p\in\mathbb{R}_{1}^{n+1}\,:\,\langle p,p\rangle =r^{2}\}$
admits different timelike, closed and conformal vector fields.
 It is a straightforward consequence of the fact that it is a
 GRW, homogeneous and isotropic. This vector fields are obtained
 fixing $p_{0}\in \mathbb{R}^{n+1}_{1}$ with $\langle p_{0},p_{0}\rangle=-1$
and taking $V_{p}=p_{0}-\frac{\langle p_{0},p\rangle}{r^{2}}p$. The different
 decompositions of the de Sitter space are
 $\mathbb{R}\times_{f} \mathbb{S} ^{n-1}(\mu)$ where
$f(t)=\frac{r}{\mu}\cosh (\frac{1}{r}t+b)$ with $r,\mu\in\mathbb{R}^{+}$ and
 $b\in\mathbb{R}$.}
\end{example}

\begin{lemma}

\label{esfera} Let $M$ be a complete Riemannian manifold with $n\geq 2$ and $%
V$ a closed and conformal vector field. Call $U$ its unitary and $A=\{p\in\
M \, : \, \lambda(p)\neq 0\}$, where $\lambda=|V|$. If the equation $%
U(U(\lambda))=-c^{2}\lambda$ holds in $A$, then $M$ is isometric to a sphere
of curvature $c^2$, $\mathbb{S}^{n}(\frac{1}{c})$.
\end{lemma}

\begin{proof}
\ It is known that $V$ has at most two zeroes, \cite{Tashiro}. Suppose $%
\nabla V=a\cdot id$. Then $\nabla a=U(a)U=U(U(\lambda))U=-c^{2}\lambda U=-c^{2}V$ and $Hess_{a}=-c^{2}a\, g$ in $M$. Using \cite{Tashiro}
(Theorem 2, III), we conclude that $M$ is a sphere of curvature $c^{2}$.
\end{proof}

The following theorem shows that the sign of the lightlike sectional
curvature is an obstruction for the existence of more than one GRW
decomposition. In addition, we show that the de Sitter space is the unique
complete Lorentzian manifold with more than one non trivial GRW
decomposition.

\begin{theorem}

\label{unicidad} Let $M$ be a Lorentzian manifold with $n\geq 3$ and $V$ a
timelike, closed and conformal vector field. If there exists another closed
and conformal vector field $W$, without zeros and linearly independent to $V$
at some point $p\in M$, then

\begin{enumerate}

\item There are degenerate planes $\Pi$ of $T_{p}M$ such that $\mathcal{K}%
_{U}(\Pi)=0$, being $U$ the unitary of $V$.

\item If $V$ is not parallel and $M$ is complete then it is isometric to a
de Sitter space.
\end{enumerate}

\end{theorem}

\begin{proof}
\ (1) Suppose $\nabla W=b\cdot id$ and $W=\alpha U+X$, where $U=\frac{V}{|V|}$ and $X\in U^{\perp}$%
. Call $-\lambda^{2}=g(V,V)$, $\sigma^{2}=g(X,X)$, $A=\{p\in M:
\sigma(p)\neq 0\}$ and define in $A$ the vector field $F=\frac{X}{\sigma}$.

Since $\nabla_{U}X=(b-U(\alpha))U$ and $0=Ug(X,U)=g(\nabla_{U}X,U)$ we have $%
b=U(\alpha)$ and $\nabla_{U}X=0$. Now, $\nabla_{X}X=\nabla_{X}(W-\alpha
U)=(b-\alpha U(\ln \lambda))X-X(\alpha)U$, and taking derivative along $X$
in $\sigma^{2}=g(X,X)$ we have $F(\sigma)=b-\alpha U(\ln\lambda)$. Thus we
get the equation $U(\alpha)=\alpha U(\ln \lambda)+F(\sigma)$. On the other
hand, $F(\alpha)=-\sigma U(\ln\lambda)$.

Taking derivative along $U$ in $\sigma^{2}=g(X,X)$ we get $U(\sigma)=0$ and
we can easily check that $[U,\lambda F]=0$. From the expression $U(\lambda
F(\alpha))=\lambda F(U(\alpha))$ we obtain the equation
\begin{equation}  \label{eqi}
F(F(\sigma))=-\sigma U(U(\ln\lambda)).
\end{equation}

Take $p\in A$ and call $\lambda(t)=\lambda(\gamma(t))$ and $%
\sigma(s)=\sigma(\zeta(s))$, where $\gamma$ is the integral curve of $U$ and
$\zeta$ of $F$ through $p$. Then $M$ is locally isometric to $%
((-\varepsilon,\varepsilon)\times L_{p},-dt^{2}+(\frac{\lambda(t)}{\lambda(0)%
})^{2}g_{\mid_{L_{p}}})$ being $\frac{\partial}{\partial t}$ identified with
$U$. Since $X$ is closed and conformal in $U^{\perp}$, $(L_{p},g_{%
\mid_{L_{p}}})$ is locally isometric to $((-\delta,\delta)\times S,ds^{2}+(%
\frac{\sigma(s)}{\sigma(0)})^{2}g_{S})$, where $\frac{\partial}{\partial s}$ is identified with $F$, see the comments after lemma
\ref{conformalvector}.

Consider the degenerate plane $\Pi=span(-U_{p}+F_p,v)$, where $v$ is an
unit vector orthogonal to $U_{p}$ and $F_{p}$. From lemma \ref%
{formulacurvnula} we obtain $\mathcal{K}_{U}(\Pi)=K^{L_p}(span(v,X)) -(\ln\lambda)^{\prime\prime}(0)$. But the curvature formulas for a warped
product gives us $K^{L_p}(span(v,X))=-\frac{\sigma^{\prime\prime}(0)}{%
\sigma(0)}$ and from equation \ref{eqi} we get $(\ln\lambda)^{\prime%
\prime}(0)=-\frac{\sigma^{\prime\prime}(0)}{\sigma(0)}$. Therefore $\mathcal{%
K}_{U}(\Pi)=0$.

(2) Suppose now that $M$ is complete and take $\gamma:\mathbb{R}\rightarrow
M $ an integral curve of $U$ with $\gamma(0)=p\in A$. Since $U(\sigma)=0$ we
have that $\gamma(t)\in A$ for all $t\in\mathbb{R}$. Then, using $[U,\lambda
F]=0$, there is a constant $k$ such that $F(F(\sigma))(\gamma(t))=-\frac{k}{%
\lambda ^2}\sigma$ for all $t$, and therefore $\lambda^{2}(\ln\lambda)^{%
\prime\prime}=k$.

If $k< 0$ the solutions of this differential equation are not positive on
the whole $\mathbb{R}$ and if $k=0$ we obtain an incomplete warped product.
If $k>0$ the solutions are $\lambda(t)=\frac{\cosh(\sqrt{k}ct+b)}{c}$ where $%
c>0 $. Since $\lambda(t)$ is not periodic, $\gamma$ can not return to $L_{p}$
and then by corollary \ref{covGRW}, $M$ is isometric to the GRW space $%
\mathbb{R}\times_{\frac{\lambda(t)}{\lambda(p)}} L_{p}$.

Now, for any $q\in L_{p}$ we have $F(F(\sigma))(q)=-\sigma (p)U(U(\ln\lambda ))(p)=-\frac{k}{\lambda^2}\sigma$ and the above lemma says that
$L_{p}$ is an euclidean sphere of curvature $\frac{k}{\lambda(p)^{2}}$. Then, comparing the warping function with that of example
\ref{seudoesfera}, $M$ is isometric to the de Sitter space of constant curvature $c^{2}k$.
\end{proof}

In the case $dim\, M=2$, $M$ complete and $V$ non parallel, the orthogonal
leaves of $V$ are trivially isometric to $(\mathbb{R},dx^{2})$ or $(\mathbb{S%
}^{1},dx^{2})$. We can obtain as in the above proof that $\lambda(t)=\frac{%
\cosh(\sqrt{k}at+b)}{b}$ and therefore $M$ is isometric to the $2$-de Sitter space or its universal covering.

The completeness hypothesis in theorem \ref{unicidad} (2) is necessary as is shown by the following example.

\begin{example}
\emph{Take $(L,g_{0})$ any Riemannian manifold and consider $M=(0,\infty)\times\mathbb{R}\times L$ with metric $
g=-dt^2+t^{2}(ds^2+e^{2s}g_{0})$. The vector fields $t\frac{\partial}{
\partial t}$ and $(t-e^{s})\frac{\partial}{\partial t}+\frac{e^{s}}{t}\frac{
\partial}{\partial s}$ are closed and conformal. Observe that the second vector field is not timelike on all the spacetime}
\end{example}

\begin{corollary}

The Friedmann Models admit an unique GRW decomposition, even locally.
\end{corollary}

\begin{proof}
\ From Friedmann equation and lemma \ref{formulacurvnula} we get $\mathcal{K}%
_{\frac{\partial }{\partial t}}(\Pi )>0$ for any degenerate plane $\Pi $.
\end{proof}

\textbf{Acknowledgments}. We thanks the referee for his valuable suggestions to improve this paper.


\begin{thebibliography}{99}
\bibitem{alias} L. J. Al\'{\i}as, A. Romero and M. S\'{a}nchez, \textit{%
Uniqueness of complete spacelike hypersurfaces of constant mean curvature in Generalized Robertson-Walker spacetimes}, Gen. Relativity
Gravitation \textbf{27} (1995), 71-84.

\bibitem{Alias2} L. J. Al\'{\i}as, A. Romero and M. S\'{a}nchez, \textit{Spacelike hypersurfaces of constant mean curvature and
Calabi-Bernstein type problems}, Tohoku Math. J. \textbf{49} (1997), 337-345.

\bibitem{Ellis} C. M. Claudel, K. S. Virbhadra and G. F. R. Ellis, \textit{%
The geometry of photon surfaces}, L. Math. Phys. \textbf{42} (2001), 818-838.

\bibitem{Esche} J. H. Eschenburg and E. Heintze, \textit{Unique
decomposition of Riemannian manifolds}, Proc. Amer. Math. Soc. \textbf{126}
(1998), 3075-3078.

\bibitem{Fernandez} L. Fern\'{a}ndez-Jambrina and R. Lazkoz, \textit{%
Classifications of cosmological milestones}, arXiv Physics gr-qc/0607073.

\bibitem{Garcia3} M. Fern\'{a}ndez-L\'{o}pez, E. Garc\'{\i}a-R\'{\i}o and D.
N. Kupeli, \textit{A characterization of cosmological time functions}, Ann.
Global Anal. Geom. \textbf{21} (2002),1-13.

\bibitem{Garcia2} E. Garc\'{\i}a-R\'{\i}o and D. N. Kupeli, \textit{%
Singularity versus splitting theorems for stably causal spacetimes}, Ann. Global Anal. Geom. \textbf{14} (1996), 301-312.

\bibitem{Harris1} S. G. Harris, \textit{A triangle comparison theorem for
Lorentz manifolds}, Indiana Univ. Math. J. \textbf{31} (1982), 289-308.

\bibitem{Harris2} S. G. Harris, \textit{A characterization of
Robertson-Walker spaces by lightlike sectional curvature}, Gen. Relativity
Gravitation \textbf{17} (1985), 493-498.

\bibitem{Koch} L. Koch-Sen, \textit{Infinitesimal null isotropy and
Robertson-Walker metrics}, J. Math. Phys. \textbf{26} (1985), 407-410.

\bibitem{Kuhnel} W. Kuhnel and H. B. Rademacher, \textit{Conformal vector
fields on pseudo-Riemann spaces}, Differential Geom. Appl. \textbf{7}
(1997),
 237-250.

\bibitem{ONeill} B. O'Neill, \textit{Semi-Riemannian geometry with
applications to relativity}, Academic Press, New York (1983).

\bibitem{Ponge} R. Ponge and H. Reckziegel, \textit{Twisted products in
pseudo-Riemannian Geometry}, Geom. Dedicata \textbf{48} (1993), 15-25.

\bibitem{Tashiro} Y. Tashiro, \textit{Complete riemannian manifolds and some
vector fields}, Trans. Amer. Math. Soc. \textbf{117} (1965), 251-275.

\bibitem{Sanchez} M. S\'{a}nchez, \textit{On the geometry of Generalized
Robertson-Walker Spacetimes: Geodesics}, Gen. Relativity Gravitation \textbf{%
30} (1998), 915-932.

\bibitem{Sanchez2} M. S\'{a}nchez and J. M. M. Senovilla, \textit{A note on the uniqueness of global static decompositions}, arXiv Physics gr-qc/0709.0305.

\end{thebibliography}
\end{document}